\documentclass[11pt]{article}
\usepackage{amsfonts}

\usepackage{amssymb}
\usepackage{latexsym}
\usepackage{amscd}
\usepackage[centertags]{amsmath}
\usepackage{xypic}

\newtheorem{thm}{Theorem}

\numberwithin{equation}{section}

\newcommand{\rmd}{\mathrm{d}}

\newcommand{\Enum}{\mathbb{E}}
\newcommand{\Pnum}{\mathbb{P}}
\newcommand{\Qnum}{\mathbb{Q}}
\newcommand{\Rnum}{\mathbb{R}}

\newcommand{\Nnum}{\mathbb{N}}

\newcommand{\lam}{\lambda}
\newcommand{\Lam}{\Lambda}
\newcommand{\vep}{\varepsilon}

\newcommand{\wt}{\widetilde}

\begin{document}
\title{Bounding the Solutions to Some SDEs via Ergodic Theory}
\author{Jian-Sheng Xie\footnote{E-mail: jiansheng.xie@gmail.com}\\
{\footnotesize School of Mathematical Sciences, Fudan
University, Shanghai 200433, P. R. China}} 
\date{}
\maketitle
\begin{abstract}
In this note we consider autonomous SDEs admitting smooth invariant measures. We present a
method in finding (almost everywhere) good bounds for $\sup \{\|X_t\|: t \in [0, T]\}$ for
strong solutions $X_{\cdot}$ to such SDEs, which in many cases are optimal bounds. In some situation
(especially in one-dimensional SDEs' cases), the discarded measure-zero set can be chosen
to be a measure-zero set of the underlying Brownian motion uniform for all initial points
$X_0=x$.
\end{abstract}

\section{Introduction}\label{sec:1}
It's well known that, for a given one-dimensional stationary Ornstein-Uhlenbeck (OU for short) process
$X=\{X_t: t \geq 0\}$£¬there exist $\lam, \sigma>0, \mu \in \Rnum$ and a standard Brownian Motion (BM for
short) $B (\cdot)$ such that $X$ has the same distribution as $\{\sigma \cdot e^{-\frac{\lam t}{2}} \cdot
B (e^{\lam t}) + \mu: t \geq 0 \}$. Therefore the law of iterated logarithm for BM (see, e.g., \cite{RM99})
leads us to the conclusion $X_t= O (\sqrt{\log t})$ almost surely. In a previous note \cite{Xie14} we have
proved the validity of such a bound for general OU processes $X=\{X_t: t \geq 0\}$ via elementary arguments
in Ito's stochastical analysis theory. In this note, we will consider the following SDE with smooth coefficients
\begin{equation}\label{eq: SDE-0}
\rmd X_t=b (X_t) \rmd t +\sigma (X_t) \rmd W_t
\end{equation}
which admits stationary strong solutions with steady distribution denoted by $\mu$; we will denote by $\Pnum$
the distribution of Wiener process $W_{\cdot}$ and $\Pnum_\mu$ the distribution of the stationary strong
solution $X_{\cdot}$. And we are interested in the growth of $\sup \{\|X_t\|: t \in [0, T]\}$ in terms of $T$,
where $\|\cdot\|$ is the Euclidean norm. We will present an ergodic theoretic method in solving such problems.
We recall that, the system (\ref{eq: SDE-0}) is called strong complete \cite{Elworthy82'}, if its solutions $X_{\cdot}$ with 
arbitrarily initial value $X_0=x$ is continuous in $t \in [0, +\infty)$ for all Wiener process orbits $W_{\cdot}$ 
in a common full standard Wiener-measure set; see, for instance, \cite{BF61} \cite{Elworthy82} \cite{Kunita90} \cite{FIZ07} 
for results relating the property of strong completeness.

Our main result may be stated as the following.
\begin{thm}\label{thm: main}
Suppose that the smooth coefficients of (\ref{eq: SDE-0}) satisfy
$$
\varlimsup_{\|x\| \to +\infty} \frac{\|b (x)\| +\|\sigma (x)\|}{\|x\|^m}<\infty
$$
for some $m \in \Nnum$. Assume the steady distribution $\mu$ is such that there exists a smooth positive
function $V$ with $\int e^{\delta V (x)} \rmd \mu<\infty$ for all $\delta \in (0, 1)$ and
$$
\lim_{\|x\| \to +\infty} \frac{V (x)}{\log \|x\|}=\infty, \quad
\varlimsup_{\|x\| \to +\infty} \frac{\|\nabla V (x)\| +\|\mathrm{Hess}_V (x)\|}{\|x\|^m}<\infty.
$$
Here $\mathrm{Hess}_V (x)$ denotes the Hessian of $V$; we always assume the monotonicity of $V (x)$ in $\|x\|$
for large $\|x\|$. Then the solution to (\ref{eq: SDE-0}) always satisfies
\begin{equation}\label{eq: statement}
\varlimsup_{t \to \infty} \frac{V (X_t)}{\log t} \leq 1
\end{equation}
$\Pnum_\mu$ almost surely.

If furthermore both $\mu$ and the transition probability semigroup of
(\ref{eq: SDE-0}) have smooth densities, then the statement (\ref{eq: statement}) holds true for all
initial values $X_0=x$ and $\Pnum$-a.e. Wiener orbits $W_{\cdot}$; in one-dimensional case with the assumption of
the strong completeness of the system, the validity of this statement can even be strengthened to be valid for all initial 
values $X_0=x$ and all Wiener process orbits $W_{\cdot}$ in a $\Pnum$-full measure set.
\end{thm}
By choosing a suitable smooth function $V (\cdot)$, it is possible to get good bounds for the growth of
$\sup \{\|X_t\|: t \in [0, T]\}$ in terms of $T$ as the examples reveal. Such result seems to be new in
literature as to our knowledge and deserves a publication somewhere.
\section{Proof of the Main Theorem}\label{sec:2}
First assume $X_{\cdot}$ to be a stationary strong solution to (\ref{eq: SDE-0}). Consider
$$
Y_t :=e^{\delta V (X_t)} \hbox{ with } \delta \in (0, \frac{1}{2}).
$$
It is clear that for suitable choice of $\wt{b} (\cdot)$ and $\wt{\sigma}$
$$
\rmd Y_t=Y_t [\wt{b} (X_t) \rmd t +\wt{\sigma} (X_t) \rmd W_t]
$$
with
$$
\varlimsup_{\|x\| \to +\infty} \frac{|\wt{b} (x)| +|\wt{\sigma} (x)|}{\|x\|^{4 m}}<\infty.
$$
This guarantees the integrability of $\wt{b} \cdot e^{\delta V}$ and $|\wt{\sigma}|^2
\cdot e^{2 \delta V}$ with respect to $\mu$. Define $M_t :=\int_0^t \wt{\sigma} (X_s) \rmd W_s$.
It is easy to see that $M_{\cdot}$ is an $L_2$-martingale with
\begin{equation}\label{eq: QV}
<M>_t :=\int_0^t |\wt{\sigma} (X_s)|^2 \cdot e^{2 \delta V (X_s)}\rmd s.
\end{equation}
By Birkhoff's ergodic theorem
$$
\lim_{t \to \infty} \frac{<M>_t}{t}=\Enum_\mu [|\wt{\sigma} (X_0)|^2 \cdot e^{2 \delta V (X_0)}]<\infty.
$$
The {\sl law of iterated logarithm} (abbr. LIL) for continuous martingale \cite{RM99} then tells us
$$
\lim_{t \to \infty} \frac{M_t}{t}=0 \hbox{ almost surely}.
$$
It is also clear that
$$
\lim_{t \to \infty} \frac{1}{t} \int_0^t Y_s \cdot \wt{b} (X_s) \rmd s =\Enum_\mu [Y_0 \cdot \wt{b}
(X_0)] \hbox{ almost surely}.
$$
Therefore $e^{\delta V (X_t)} \leq a \cdot t +b$ for all $t \geq 0$ almost surely, where $a$ is a positive
constant and $b$ is a measurable function of $X_{\cdot}$ independent of $t$. Hence
\begin{equation}\label{eq: statement'}
\varlimsup_{t \to \infty} \frac{V (X_t)}{\log t} \leq 2
\end{equation}
$\Pnum_\mu$ almost surely.

Now we are going to lower the bound 2 in the right hand side of (\ref{eq: statement'}) down into 1 as (\ref{eq: statement})
says. Assume that we have already proved
$$
\varlimsup_{t \to \infty} \frac{V (X_t)}{\log t} \leq \beta
$$
for some constant $\beta$; therefore for any $\vep \in (0, 1)$ there is $C$ such that
$$
e^{V (X_t)} \leq C (t^{\beta+\vep} +1), \quad \forall t \geq 0.
$$
Now fix a number $\delta \in (1/2, 1)$ arbitrarily and define $Y_t, M_t$ as above.
Then using the above arguments once again, we find
\begin{eqnarray*}
\frac{<M>_t}{t} &=& \frac{1}{t} \int_0^t \wt{\sigma} (X_s)^2 e^{(2 \delta +\vep -1) V (X_s)} \cdot e^{ (1-\vep) V (X_s)} \rmd s\\
&\leq& \frac{1}{t} \int_0^t e^{(1-\vep) V (X_s)} \rmd s \cdot [C (t^{\beta+2 \vep} +1)]^{2 \delta+2 \vep -1}
\end{eqnarray*}
for sufficiently large $t$. In view of LIL for continuous martingale \cite{RM99} and Birkhoff's ergodic theorem, this implies
$$
\lim_{t \to +\infty} \frac{M_t}{t^{\frac{1}{2}+(\delta+2 \vep -\frac{1}{2}) \cdot (\beta+2\vep)}}=0.
$$
On the other hand, we still have $\wt{b} \cdot e^{\delta V} \in L_1 (\mu)$. Therefore we have
$$
\lim_{t \to +\infty} \frac{e^{\delta V (X_t)}}{t^{\gamma}}=0
$$
for $\gamma :=\max \{ 1+\vep, \frac{1}{2}+(\delta+2 \vep -\frac{1}{2}) \cdot (\beta+2\vep) \}$. This proves
$$
\varlimsup_{t \to +\infty} \frac{V (X_t)}{\log t} \leq \frac{\gamma}{\delta}.
$$
Letting $\delta \to 1$ and then $\vep \to 0$, we obtain
$$
\varlimsup_{t \to +\infty} \frac{V (X_t)}{\log t} \leq \beta^\prime :=\max\{1, \frac{1+\beta}{2}\}=\frac{1+\beta}{2}
$$
with initial $\beta=2$. This machinery leads us finally to (\ref{eq: statement}).

Qian and Zhang's argument \cite[page 1637]{QZ05} tells us that, when $\mu$ and the transition probability semigroup
of $X_{\cdot}$ have densities, (\ref{eq: statement}) holds true for all $X_0=x$ and all
Wiener orbits $W_{\cdot} \in \Lam_x$ with $\Pnum (\Lam_x)=1$. In one dimensional case, we can say more.
Let
$$
\Lam :=\bigcap_{r \in \Qnum} \Lam_r,
$$
so $\Pnum (\Lam)=1$. For any two solutions $X^x_{\cdot}, X^y_{\cdot}$ to (\ref{eq: SDE-0}) with initial values $X_0^x=x,
X_0^y=y, x \neq y \in \Rnum$, write $Z_t:= X_t^x -X_t^y$. It is easy to see that
$$
\rmd Z_t=Z_t \cdot [\wt{b} (X_t^x, X_t^y) \rmd t +\wt{\sigma} (X_t^x, X_t^y) \rmd W_t]
$$
for some smooth functions $\wt{b} (x, y), \wt{\sigma} (x, y)$. Then one clearly has
$$
Z_t=Z_0 \cdot \exp \Bigl( \int_0^t[\wt{b} (X_s^x, X_s^y) -\frac{1}{2} \wt{\sigma} (X_s^x, X_s^y)^2] \rmd s
+\int_0^t \wt{\sigma} (X_s^x, X_s^y) \rmd W_s \Bigr),
$$
which implies that
\begin{equation}
X_t^x \leq X_t^y \hbox{ for all } t \geq 0 \hbox{ if } x<y.
\end{equation}
If $\varphi$ is a continuous monotonic function in $L_1 (\mu)$, then it is easy to see
that Birkhoff's ergodic theorem holds
$$
\lim_{t \to \infty} \frac{1}{t} \int_0^t \varphi (X_s) \rmd s =\int \varphi (x) \rmd \mu
$$
for all $x \in \Rnum$ and all $W_{\cdot} \in \Lam$. The same argument applies to the limit
$$
\varlimsup \frac{V (X_t)}{\log t} \leq 1
$$
yielding its validity for all $x \in \Rnum$ and all $W_{\cdot} \in \Lam$, if $V (x)$ is increasing
in $|x|$ for large $|x|$.
\section{Examples and Discussions}\label{sec:3}
The first example is the standard OU process $\{X_t :=e^{-t} B (e^{2 t}): t \geq 0 \}$, where
$B (\cdot)$ is a standard BM. By LIL for BM, we have
$$
\varliminf_{t \to \infty} \frac{X_t}{\sqrt{2 \log t}}=-1, \varlimsup_{t \to \infty} \frac{X_t}{\sqrt{2 \log t}}=1.
$$
While our argument in this note yields
$$
\varlimsup_{t \to \infty} \frac{|X_t|}{\sqrt{2 \log t}} \leq 1.
$$
Hence this example indicates that our method can give optimal bounds in some cases (hopefully always so).

The second example is the following one dimensional SDE:
$$
\rmd X_t =-U^\prime (X_t) \rmd t +\sqrt{2 \vep} \rmd W_t
$$
with
$$
\int e^{-U (x)/\vep} \rmd x<\infty.
$$
Here $U (x)$ is a polynomial with leading term being $c \, x^{2 p}$ for some $c>0, p \geq 1$. Then our argument gives
$$
\varlimsup_{t \to \infty} \frac{|X_t|}{(\log t)^{1/2p}} \leq (\frac{\vep}{c})^{1/2p}
$$
for all $X_0=x$ and all BM orbits $W_{\cdot} \in \Lam$ with $\Pnum (\Lam)=1$, since the strong completeness of the model 
is guaranteed by \cite{FIZ07}.

We would like to give some discussions. As is well known, Birkhoff's ergodic theorem is an extension of Kolmogrov's
{\sl strong law of larg numbers} (abbr. SLLN); It also holds for stationary processes with continuous-time parameter
under suitable $L_1$-integrabilty condition. In probability theory, when the $L_1$-integrabilty condition is replaced
by $L_p$-integrability condition (with $p \in (0, 2)$), Marcinkiewicz-Zygmund's SLLN would take place of Kolmogrov's
SLLN for i.i.d. random variables sequence. It is easy to see that the following result holds, which is a generalization 
of one part of Marcinkiewicz-Zygmund's SLLN.
\begin{thm}
Let $\{X_n: n \geq 0\}$ be an stationary process with $\Enum |X_0|^p<\infty$ for some $p \in (0, 1)$. Then
$$
\lim_{n \to +\infty} \frac{1}{n^{1/p}} \sum_{k=0}^{n-1} X_k=0 \quad \hbox{ almost surely.}
$$
\end{thm}
The counterpart of $p \in (1,2)$ as the above theorem to Marcinkiewicz-Zygmund's SLLN seems still unknown. Also, it is
interesting to ask the validity of the continuous-time counterpart of the above theorem. It seems to us that, a proper statement
might be as the following: Let $\{X_t: t \geq 0\}$ be an stationary process with $\Enum |X_0|^p<\infty$ for some $p \in (0, 1)$.
Assume the continuity of $X_t$ in $t$. Then for all $\vep>0$
$$
\lim_{T \to +\infty} \frac{1}{T^{1/p+\vep}} \int_0^T X_s \rmd s=0 \quad \hbox{ almost surely.}
$$
It is the deficiency of such a result that forces us to find the machinery mentioned in the second paragraph
of Sect. \ref{sec:2}.

\noindent{\sl \textbf{Acknowledgements} \quad} {The author thanks Prof. Jiangang Ying for helpful discussions. He
is also grateful to Prof. Xi Chen for telling him some literatures on the strong completeness of SDEs.
This work is partially supported by NSFC (No. 10701026 and No. 11271077) and the Laboratory of Mathematics for
Nonlinear Science, Fudan University.}






\begin{thebibliography}{00}


\bibitem{BF61} Blagove\v{s}\v{c}enski\u{\i}, Ju. N.; Fre\u{\i}dlin, M. I.:
{\em Some properties of diffusion processes depending on a parameter. (Russian)}
Dokl. Akad. Nauk SSSR \textbf{138} (1961), pp. 508¨C-511.

\bibitem{Elworthy82} Elworthy, K. D.:
{\em Stochastic flows and the C0-diffusion property.}
Stochastics \textbf{6} (1981/82), no. 3-4, pp. 233¨C-238.

\bibitem{Elworthy82'} Elworthy, K. D.:
{\em Stochastic differential equation on manifolds.} 
London Mathematical Society Lecture Note Series, 70. Cambridge University Press, Cambridge-New York, 1982. xiii+326 pp.


\bibitem{FIZ07} Fang, Shizan; Imkeller, Peter; Zhang, Tusheng:
{\em Global flows for stochastic differential equations without global Lipschitz conditions.}
Ann. Probab. \textbf{35} (2007), no. 1, pp. 180-¨C205.

\bibitem{Kunita90} Kunita, Hiroshi:
{\em Stochastic flows and stochastic differential equations.}
Cambridge Studies in Advanced Mathematics, \textbf{24}. Cambridge University Press,
Cambridge, 1990. xiv+346 pp. ISBN: 0-521-35050-6

\bibitem{QZ05} Qian, M.; Zhang, Fu-Xi:
{\em Non-equilibrium of a general stochastic system of coupled oscillators:
entropy production rate and rotation numbers}, Ergod. Theory \& Dynam. Sys. 25,
no. 5 (2005), pp. 1633--1641.

\bibitem{RM99} Revuz, Daniel; Yor, Marc:
{\it Continuous martingales and Brownian motion}.
Third edition. Grundlehren der Mathematischen Wissenschaften
[Fundamental Principles of Mathematical Sciences], 293.
Springer-Verlag, Berlin, 1999. xiv+602 pp. ISBN: 3-540-64325-7

\bibitem{Xie14} Xie, Jian-Sheng:
{\em Bounding Ornstein-Uhlenbeck Processes and Alikes}, preprint.

\bibitem{Zhang10} Zhang, Xicheng,
{\em Stochastic Homeomorphism Flows of SDEs with Singular Drifts and Sobolev Diffusion Coefficients.}
arXiv:1010.3403.
\end{thebibliography}



\end{document}